\documentclass[12pt]{article}
\usepackage{amsthm,amsmath}
\usepackage{amscd}
 \usepackage{amssymb}
 \usepackage{hyperref}
 \usepackage[all]{xypic}
 \usepackage{mathtools}
 \usepackage[utf8]{inputenc}
  \newtheorem{lemma}{Lemma}[section]
 \newtheorem{corollary}[lemma]{Corollary}
 
 \newtheorem{theorem}[lemma]{Theorem}
 \newtheorem{proposition}[lemma]{Proposition}
 
 \newtheorem{remark}[lemma]{Remark}

\usepackage[utf8]{inputenc} %unicode support
\def\includegraphics{}

\usepackage{tikz,xcolor,hyperref}

\definecolor{lime}{HTML}{A6CE39}
\DeclareRobustCommand{\orcidicon}{
	\begin{tikzpicture}
	\draw[lime, fill=lime] (0,0)
	circle[radius=0.16]
	node[white]{{\fontfamily{qag}\selectfont \tiny \.{I}D}};
	\end{tikzpicture}
	\hspace{-2mm}
}
\foreach \x in {A, ..., Z}{%
	\expandafter\xdef\csname orcid\x\endcsname{\noexpand\href{https://orcid.org/\csname orcidauthor\x\endcsname}{\noexpand\orcidicon}}
}

\usepackage{amscd}
\usepackage{amsmath,amssymb,mathrsfs}
\usepackage{hyperref}
\usepackage[all]{xypic}
\usepackage{times}
\usepackage{inputenc}
\usepackage{titlesec}

\allowdisplaybreaks[1]
\allowdisplaybreaks[4]

 \usepackage{amsthm}

\advance\headheight1.15pt

\newenvironment{proof of Theorem 3.1}{{\noindent\it {\bf Proof of Theorem 3.1}}\quad}{\hfill $\square$\par}
\newenvironment{proof of Theorem 2.9}{{\noindent\it {\bf Proof of Theorem 2.9}}\quad}{\hfill $\square$\par}
\newenvironment{proof of Theorem 4.1}{{\noindent\it {\bf Proof of Theorem 4.1}}\quad}{\hfill $\square$\par}
\newenvironment{proof of Theorem 1.4}{{\noindent\it {\bf Proof of Theorem 1.4}}\quad}{\hfill $\square$\par}
\newenvironment{proof of Theorem 1.5}{{\noindent\it {\bf Proof of Theorem 1.5}}\quad}{\hfill $\square$\par}
\newcommand{\I}{\int_{0}^{1}}      
\newcommand{\sk}{\sum_{k=0}^{\infty}}   \newcommand{\hp}{H^{p}}
\newcommand{\sn}{\sum_{n=0}^{\infty}}  
\newcommand{\hd}{H(\mathbb{D})}            
\newcommand{\dd}{\mathbb{D}}                     \newcommand{\pd}{D^{p}_{p-1}}
 \newcommand{\ssp}{S^{p}}   \newcommand{\sq}{S^{q}}
  \newcommand{\h}{\mathcal {H}} \newcommand{\xp}{X_{p}}
\newcommand{\hg}{\mathcal {H}_{g}} \newcommand{\ce}{\mathcal{C}_\eta}   \newcommand{\tto}{\rightarrow}
\newcommand{\comment}[1]{}

\begin{document}

\baselineskip=8pt
%\title{A Generalization of Carath\'{e}odory Theorem Using Best Proximity Theory}
\title{\fontsize{15}{0}\selectfont  Operators of  Hilbert  type  acting   on some spaces of analytic functions}
\author{\fontsize{11}{0}\selectfont
	Pengcheng Tang$^{*,a}$
	%\hspace{-1.5mm} \orcidA{}
	\\
	\fontsize{10}{0}\it{} %\fontsize{10}{0}\it{Email:@gmail.com}
	\\ \fontsize{10}{0}\it{  $^{a}$School of Mathematics and Statistics, Hunan University of Science and Technology,}
	\\ \fontsize{10}{0}\it{ Xiangtan, Hunan 411201, China}
}

\date{}
%\affil{\fontsize{12}{0}\selectfont aaaa}
%\affiliation
%\emailaddress{abhikdigar@gmail.com}
%\thanks{E-mail : abhikdigar@gmail.com}
\maketitle
\thispagestyle{empty}
\begin{center}

\textbf{\underline{ABSTRACT}}
\end{center}
\ \ \ \ \ Let  $\hd$ be the space of all analytic functions in the unit disc $\dd$. For $g\in \hd$, the  generalized
Hilbert operator  $\hg$  is defined by
$$\hg(f)(z)=\int_{0}^{1}f(t)g'(tz)dt, \ \ z\in \dd, f\in \hd.$$

In this paper, we study the operator  $\hg$  acting on some  spaces of analytic functions in $\dd$. Specifically,  we give a complete characterization of those $g\in \hd$ for which  the operator $\hg$ is bounded (resp. compact)   from  the Dirichlet  space $\mathcal{D}^{2}_{\alpha}$ to $\mathcal{D}^{2}_{\beta}$ for all possible indicators $\alpha,\beta \in \mathbb{R}$. We also study the action of the  operator $\hg$ on the space of bounded analytic functions $H^{\infty}$, which generalizes the known results for the classical Hilbert operator $\mathcal {H}$ acting on $H^{\infty}$. In particular, we consider the boundedness of the operator  $\hg$
with a symbol of non-negative Taylor coefficients, acting on logarithmic Bloch spaces and on Korenblum spaces. This work generalizes the corresponding results for the classical Hilbert operator.  
%Additionally,  we  use our results concerning the operator  $\hg$  to study certain Ces\`{a}ro-type operators   $\mathcal C_{(\eta )}$ acting on the Wiener algebra and on derivative Hardy spaces. 

\iffalse
 For $0<p<\infty$, we let  $S^{p}$
be the space of functions $f\in \hd$ such that $f'$  belongs to the Hardy spaces $H^{p}$.
Let $(\eta )=\{ \eta _n\}_{n=0}^\infty $  be  a sequence of
complex numbers and $f(z)=\sum_{n=0}^\infty
a_nz^n\in \hd$,  the Ces\`{a}ro type operator  $\mathcal C_{(\eta )}$ is formally defined by
$$\mathcal C_{(\eta )}(f)=\mathcal C_{\{
	\eta_n\}}(f)(z)=\sum_{n=0}^\infty \eta _n\left (\sum
_{k=0}^na_k\right )z^n.$$ The operator $\mathcal C_{(\eta )}$ is a
natural generalization of the Ces\`{a}ro operator.

In this paper, we provide a characterization of those  $\{\eta_{n}\}$ for which  $\mathcal C_{(\eta )}$
is bounded (resp. compact) between derivative Hardy spaces. This partially answers the question posed by Lin  and Xie in [J. Funct. Anal. 288(6) (2025) Paper No. 110813].
Alternatively, we  use our results concerning the operators  $\mathcal C_{(\eta )}$ to study certain Hilbert-type operators  $\mathcal {H}_{g}$ acting
between the derivative Hardy spaces.  We also study the range of $\hg$ acting on the space of bounded analytic functions  $H^{\infty}$. This is a generalization of the known results for the classic Hilbert operators  $\mathcal {H}$ on  $H^{\infty}$.
\fi

\let\thefootnote\relax\footnote{$^*$Corresponding Author}
\let\thefootnote\relax\footnote{ Pengcheng Tang: www.tang-tpc.com@foxmail.com}
\vspace{1cm}
% Fixed point theory has applications in almost all branches of mathematics.
% The necessary condition, for a map $T : A \rightarrow B$ to admit a fixed point,
% is $A \cap B \neq \emptyset.$ Let $A$ and $B$ be subsets of a metric space $(X,d)$.
% For a map, $T : A \cup B \rightarrow A \cup B$, we aim to find $x \in A$, $y \in B$
% such that $Tx=x,~Ty=y$ and $d(x, y) = \dist (A, B)$. In \cite{ruc},
% Eldred {\it{et al.}} introduced a class of mappings called {\it relatively u-continuous
% mappings} and proved the  existence of a point $x \in A$ such that $d(x,Tx)=\dist(A,B)$
% satisfying $TA \subset B$ and $TB \subset A$ (such a point is said to be a best proximity
%% point for $T$). In this talk, we give sufficient conditions for the existence of fixed
% points $x$ and $y$ in $A$ and $B$ respectively with $d(x,y)=\dist (A,B)$ for a relatively
% u-continuous map $T$ satisfying $TA \subset A$ and $TB \subset B$. As a consequence, we
% get the main result of \cite{ruc}. We also prove that every relatively u-continuous map
% $T$ is continuous in setting of a strictly convex Banach space. Finally, we discuss the
% existence of geranalized equilibrim pairs for generalized two person games.
%\begin{flushleft}
%
%\textbf{Keywords : } Herz spaces; Lorentz spaces; Banach function spaces; Sublinear operators.
%\end{flushleft}

\section{Introduction} \label{Sec:Intro}

\ \ \  \ \ \ Let $\mathbb{D}=\{z\in \mathbb{C}:\vert z\vert <1\}$ denote the open unit disk of the complex plane $\mathbb{C}$ and $H(\mathbb{D})$ denote the space of all
analytic functions in $\mathbb{D}$.

The Bloch type space $\mathcal {B}^{\alpha}$ consists of those functions $f\in  H(\mathbb{D}) $ for which
$$
\vert \vert f\vert \vert _{\mathcal {B}^{\alpha}}=\vert f(0)\vert +\sup_{z\in \mathbb{D}}(1-\vert z\vert ^{2})^{\alpha}\vert f'(z)\vert <\infty.
$$
The space $\mathcal {B}^{1}$ is just the classic Bloch space  and denoted simply as   $\mathcal {B}$.

%For $0<\alpha<\infty$, the Korenblum space $H_\alpha^{\infty}$ is the space of all functions $f \in H(\mathbb{D})$ for which
%$$
%\|f\|_{H_{\alpha}^{\infty}}=\sup _{z \in \mathbb{D}}(1-|z|)^\alpha|f(z)|<\infty \text {. }
%$$

Let  $0<p\leq\infty$, the classical Hardy space $H^p$  consists of  those functions  $g\in H(\dd)$ for which
$$
||g||_{H^{p}}=\sup_{0\leq r<1} M_p(r, g)<\infty,
$$
where
$$
M_p(r, g)= \left(\frac{1}{2\pi}\int_0^{2\pi}|g(re^{i\theta})|^p d\theta \right)^{1/p}, \ 0<p<\infty,
$$
$$M_{\infty}(r, g)=\sup_{|z|=r}|g(z)|.$$

For $0<p\leq\infty$, the derivative Hardy space  $S^{p}$ consists of  those functions  $g\in H(\dd)$ for which
$$||g||^{p}_{S^{p}}=|g(0)|^{p}+||g'||^{p}_{H^p}<\infty.$$
The space $S^{2}$
is a  Hilbert space and if $g(z)=\sn b_{n}z^{n}\in \hd$, then
$$||g||^{2}_{S^{2}}=|b_{0}|^{2}+\sum_{n=1}^{\infty}n^{2}|b_{n}|^{2}.$$

%For $0<p<\infty$ and $\alpha>-1$,  the weighted Bergman space $A^{p}_{\alpha}$ consists of those
%$f\in \hd$ such that
%$$||f||_{A^{p}_{\alpha}}=\left((\alpha+1)\Id |f(z)|^{p}(1-|z|^{2})^{\alpha}dA(z)\right)^{\frac{1}{p}}<\infty$$

%For $0<p<\infty$ and $\alpha>-1$, the weighted Dirichlet space  $\mathcal {D}^{p}_{\alpha}$
%consists of those
% $f\in \hd$ such that
%%$$||f||^{p}_{\mathcal {D}^{p}_{\alpha}}=|f(0)|^{p}+||f'||^{p}_{A^{p}_{\alpha}}<\infty.$$
%If $p<\alpha+1$, then  $\mathcal {D}^{p}_{\alpha}=A^{p}_{\alpha-p}$. If $p>\alpha+2$, then $\mathcal {D}^{p}_{\alpha}\subset H^{\infty}$. The space $\mathcal {D}^{2}_{0}$
% is the classical Dirichlet space  $\mathcal {D}$ and we shall write  $||f||_{\mathcal {D}}$ for  $||f||_{\mathcal {D}^{2}_{0}}$.
%The spaces $\mathcal {D}^{2}_{\alpha}(\alpha>-1)$  are
% Hilbert spaces and  if  $f(z)=\sn a_{n}z^{n}\in \hd$, then
%$$||f||^{2}_{\mathcal {D}^{2}_{\alpha}}\asymp|a_{0}|^{2}+\sum_{n=1}^{\infty}n^{1-\alpha}|a_{n}|^{2}.$$

For $\alpha\in\mathbb{R}$, we  use  $\mathcal {D}^{2}_{\alpha}$ to denote the space of functions $g(z)=\sn b_{n}z^{n}\in\hd$ such that
$$||g||^{2}_{\mathcal {D}^{2}_{\alpha}}=|b_{0}|^{2}+\sum_{n=1}^{\infty}n^{1-\alpha}|b_{n}|^{2}<\infty.$$
In particular, $\mathcal {D}^{2}_{1}=H^{2}$ and $\mathcal {D}^{2}_{-1}=S^{2}$. The space $\mathcal {D}^{2}_{0}$
is the classical Dirichlet space  $\mathcal {D}$ and we shall write  $||g||_{\mathcal {D}}$ for  $||g||_{\mathcal {D}^{2}_{0}}$.

The analytic Wiener algebra $\mathcal {W}$ consists of those function $f(z)=\sn a_{n}z^{n}\in \hd$ for which
$$||f||_{\mathcal {W}}=\sn |a_{n}|<\infty.$$

Let  $1\leq p<\infty$ and $0<\alpha\leq 1$, the mean Lipschitz space $\Lambda^p_\alpha$ consists of those functions $f\in H(\dd)$ having a non-tangential limit  almost everywhere such that $\omega_p(t, f)=O(t^\alpha)$ as $t\to 0$. Here $\omega_p(\cdot, f)$ is the integral modulus of continuity of order $p$ of the function $f(e^{i\theta})$. It is  known (see \cite{b1}) that $\Lambda^p_\alpha$ is a subset of $H^p$ and
$$\Lambda^p_\alpha=\left\{f\in H(\dd):M_p(r, f')=O\left(\frac{1}{(1-r^{2})^{1-\alpha}}\right), \ \ \mbox{as}\ r\rightarrow 1\right\}.$$

The space of those  $f\in \hd$ such that
$$M_p(r, f')=o\left(\frac{1}{(1-r^{2})^{1-\alpha}}\right) ,\ \ \ \mbox{as} \ r \tto 1^-, $$
is denoted by $\lambda^{p}_{\alpha}$.

The Hilbert matrix  $\mathcal {H}$ is an infinite matrix  whose entries are $a_{n,k}=(n+k+1)^{-1}$.
\[
\mathcal {H}=\begin{pmatrix}
1& \frac{1}{2} &\frac{1}{3}& \frac{1}{4}& \cdots\\
\frac{1}{2}& \frac{1}{3} & \frac{1}{4}&  \frac{1}{5}& \cdots\\
\frac{1}{3}& \frac{1}{4} &  \frac{1}{5}&  \frac{1}{6}& \cdots\\
\vdots& \vdots& \vdots& \vdots& \ddots
\end{pmatrix}.
\]
The matrix  $\h$ induces an operator on $\hd$ by its action on the Taylor
coefficients:
$$a_{n}\rightarrow\displaystyle{\sum_{k=0}^{\infty}}\mu_{n,k}a_{k},\  n\in \mathbb{N}\cup \{0\}.$$

The  Hilbert operator $\mathcal {H}$ defined on $\hd$ as follows:
If $f(z)=\displaystyle{\sum_{n=0}^{\infty}}a_{n}z^{n}\in \hd$,  then
\begin{equation*}\label{eq1}
\mathcal {H}(f)(z)=\sum_{n=0}^{\infty}\left(\sum_{k=0}^{\infty}\frac{a_{k}}{n+k+1}\right)z^{n}, \ \ z\in \dd,
\end{equation*}
whenever the right hand side makes sense and defines an analytic function in $\dd$.

The study of the Hilbert operator $\mathcal{H}$ on  analytic function spaces was initiated by Diamantopoulos and Siskakis in \cite{dia1}. They proved that  $\mathcal {H}:H^{p}\rightarrow H^{p}$ is bounded for $1<p<\infty$ and  $\h$ is not bounded on $H^{1}$.  Subsequently, Diamantopoulus \cite{dia2} also considered the boundedness of $\mathcal {H}$  on the Bergman spaces $A^{p}$. He proved that $\h:A^{p}\rightarrow A^{p}$ is bounded for $2<p<\infty$ and $\h$ is not bounded on $A^{2}$.
Jevti\'{c} and Karapretovi\'{c} \cite{jk}  investigated  the boundedness of $\h$ on mixed-norm spaces.  The reader is referred to \cite{bh,bell, bok,dos} for more about  Hilbert operator $\h$ on spaces of analytic functions.

If $f\in H^{1}$, then the Fej\'{e}r-Riesz inequality (see \cite[Page 46]{b1}) shows that
$\int_{0}^{1}|f(t)|dt<\infty$. This implies  that $\h (f)$ has the integral form
$$\h(f)(z)=\sn \left(\int_{0}^{1}t^{n}f(t)dt\right)z^{n}=\int_{0}^{1}\frac{f(t)}{1-tz}dt, \ \ z\in \dd.$$
Alternatively, this can be regarded as
$$\hg(f)(z)=\int_{0}^{1}f(t)g'(tz)dt \ \ \mbox{with}\ \ g(z)=\log\frac{1}{1-z}.$$

Inspired by this integral representation, the following generalized Hilbert operator $\hg$ is considered by Galanopoulos et al. \cite{ann}.
For any given  $g\in \hd$, the generalized Hilbert operator $\hg$ is defined by
\begin{equation}\label{h1}
\hg(f)(z)=\int_{0}^{1}f(t)g'(tz)dt.
\end{equation}

For any $g\in \hd$, if $f\in H^{1}$, then the integral  in (\ref{h1}) converges absolutely, and hence $\hg(f)$ is a well defined analytic function in $\dd$. Thus, if $f(z)=\sn a_{n}z^{n}\in H^{1}$
and $g(z)=\sn b_{n}z^{n}\in \hd$, then  $\hg(f)$ has the following expression in terms of Taylor coefficients:
\begin{align}\label{h}
\hg(f)(z)&=\sn \left((n+1)b_{n+1}\int_{0}^{1}t^{n}f(t)dt\right)z^{n} \nonumber\\
& = \sn \left((n+1)b_{n+1}\sk \frac{a_{k}}{n+k+1}\right)z^{n}.
\end{align}

In  \cite{ann}, the authors  studied the boundedness of $\hg$ on the Hardy spaces $H^{p}$, the Bergman spaces  $A^{p}_{\alpha}$ and on the Dirichlet type spaces $\mathcal{D}_{\alpha}^{p}$.
Pel\'{a}ez and  R\"{a}tty\"{a} \cite{pj} also investigated the generalized Hilbert operator  $\hg$
acting on the weighted Bergman space  $A^{p}_{\omega}$, where $\omega$  belongs to a specific class of regular radial weight functions. The mean Lipschitz space plays a foundational role in these works.

Recently,  Galanopoulos and Girela  used the Suchur test in  \cite{rra} to establish the boundedness and compactness of $\hg$ on  the Dirichlet space $\mathcal{D}$. In this paper, we give a complete characterization of those $g\in \hd$ for which the operator  $\hg$ is bounded (compact) from $\mathcal{D}^{2}_{\alpha}$  to $\mathcal{D}^{2}_{\beta}$  for all possible indicators $\alpha,\beta \in \mathbb{R}$. This will be the main results of the Sect. \ref{se2}. In Sect. \ref{se3}, we will be mainly devoted to study  the range of $\hg$ acting on  space of bounded analytic functions $H^{\infty}$.  In Sect. \ref{se4}, we consider the boundedness of the operator  $\hg$
with a symbol of non-negative coefficients, acting on logarithmic Bloch spaces and on Korenblum spaces. These works generalizes the corresponding results for the classical Hilbert operator.

Throughout the paper, the letter $C$ will denote an absolute constant whose value depends on the parameters
indicated in the parenthesis, and may change from one occurrence to another. We will use
the notation $``P\lesssim Q"$ if there exists a constant $C=C(\cdot) $ such that $`` P \leq CQ"$, and $`` P \gtrsim Q"$ is
understood in an analogous manner. In particular, if  $``P\lesssim Q"$  and $ ``P \gtrsim Q"$ , then we will write $``P\asymp Q"$.

 \section{ The operators $\hg$ acting between Dirichlet type spaces}\label{se2}

%  We begin with the following theorem.

\ \ \ \ \ \ If $\gamma>\alpha$, then $\mathcal{D}^{2}_{\alpha}\subset \mathcal{D}^{2}_{\gamma} $. This shows  that $\mathcal{D}^{2}_{2}\subset \mathcal{D}^{2}_{\gamma} $ for every $\gamma\geq 2$.
The function $f$ given by 
$$f(z)=\sn \frac{z^{n}}{\log(n+1)},\ \ \ z\in \dd,$$
belongs to  $\mathcal{D}^{2}_{2}$. However, we can  easily see that if $g$ is the monomial $g(z)=z^{n}(n\in \mathbb{N})$,
then  $\hg(f)$ is not well defined. Hence, $\hg(f)$ is not well defined on $\mathcal{D}^{2}_{\gamma}$ for all $\gamma\geq 2$.

If $\alpha>-1$, it is known that $f\in \mathcal{D}^{2}_{\alpha}$ if and only if 
\begin{equation}\label{f}
\int_{\dd}|f'(z)|^{2}(1-|z|^{2})^{\alpha}dA(z)<\infty.
\end{equation}
Below, we shall use  (\ref{f}) to show that  $\hg$ is  well defined on  $\mathcal{D}^{2}_{\alpha}$ for $0<\alpha<2$.
\begin{proposition}
	Let  $g(z)=\sn b_{n}z^{n}\in \hd$ and let $0<\alpha<2$. Then    the integral  $\hg(f)$ is a well defined analytic function in $\dd$ for every  $f \in  \mathcal{D}^{2}_{\alpha}$ and (\ref{h}) holds.
\end{proposition}
\begin{proof}
	It is suffices to prove that $\int_{0}^{1}|f(t)|dt<\infty$. For $f \in  \mathcal{D}^{2}_{\alpha}$,  using (\ref{f}) and the well known pointwise estimate, we have
	$$|f'(z)|\lesssim \frac{1}{(1-|z|)^{\frac{\alpha+2}{2}}}.$$
This means that 
\begin{align*}
|f(z)| &\leq |f(0)|+\int_{0}^{1}|zf'(tz)|dt \\
& \lesssim 1+|z|\int_{0}^{1}\frac{1}{(1-t|z|)^{1+\frac{\alpha}{2}}}dt\\
& \lesssim 1 +\frac{1}{(1-|z|)^{\frac{\alpha}{2}}}.
\end{align*} 	
	Since $0<\frac{\alpha}{2}<1$,this implies that 
	$$\int_{0}^{1}|f(t)|dt\lesssim 1+ \int_{0}^{1}(1-t)^{-\frac{\alpha}{2}}dt\lesssim 1.$$
The proof is complete.\end{proof}

For $0<\alpha<2$, by a theorem of    Diamantopoulos   \cite[Theorem 1.2]{dia3}, we know that the Hilbert operator $\mathcal{H}$ is bounded on $\mathcal{D}^{2}_{\alpha}$. In the next, we will apply this result to characterize those $g\in \hd$ such that  $\hg$ is bounded (resp. compact) from  $\mathcal{D}^{2}_{\alpha}$ to  $\mathcal{D}^{2}_{\beta}$ for $0<\alpha<2$ and $\beta \in R$. 
\begin{theorem}\label{th2}
	Let  $g(z)=\sn b_{n}z^{n}\in \hd$. If $0<\alpha<2$ and $\beta \in \mathbb{R}$,
	then the following conditions are equivalent:
	
	(1) The operator $\hg:\mathcal {D}^{2}_{\alpha} \tto  \mathcal {D}^{2}_{\beta}$ is bounded.
	
	(2) $\sum_{n=2^{N}-1}^{2^{N+1}}|b_{n+1}|^{2}=O\left( 2^{N(\beta-\alpha-1)}\right) $.
\end{theorem}

\begin{proof}
	$(1)\Rightarrow (2).$  Assume that $\hg:\mathcal {D}^{2}_{\alpha} \tto  \mathcal {D}^{2}_{\beta}$ is bounded. For $\frac{1}{2}<a<1$, let 
	$$f_{a}(z)=(1-a)^{\frac{\alpha}{2}}\sum_{n=0}^{\infty}(n+1)^{\alpha-1}a^{n}z^{n},\ \ \ z\in \dd.$$
	Then $||f_{a}||_{\mathcal {D}^{2}_{\alpha}}\asymp1$ for all $\frac{1}{2}<a<1$. For $N\geq 2$,  we have 
	\begin{align*}
	||\hg(f_{a})||^{2}_{\mathcal {D}^{2}_{\beta}} &\gtrsim (1-a)^{\alpha}\sum_{n=0}^{\infty}(n+1)^{3-\beta}|b_{n+1}|^{2} \left( \sum_{k=0}^{\infty}\frac{(k+1)^{\alpha-1}a^{k}}{n+k+1}\right) ^{2}\\
	& \gtrsim (1-a)^{\alpha}\sum_{n=2^N-1}^{2^{N+1}}(n+1)^{3-\beta}|b_{n+1}|^{2}\left( \sum_{k=0}^{2^{N+1}}\frac{(k+1)^{\alpha-1}a^{k}}{n+k+1}\right) ^{2}\\
	& \gtrsim (1-a)^{\alpha}\sum_{n=2^N-1}^{2^{N+1}}2^{(1-\beta)N}|b_{n+1}|^{2}\left( \sum_{k=0}^{2^{N+1}}(k+1)^{\alpha-1}\right) ^{2}a^{2^{N+1}}\\
	& \gtrsim (1-a)^{\alpha}2^{N(2\alpha-\beta-1)}a^{2^{N+1}}\sum_{n=2^N-1}^{2^{N+1}}|b_{n+1}|^{2}.
	\end{align*}
	Taking $a=1-\frac{1}{2^{N}}$, we have that 
	$$||\hg(f_{a})||^{2}_{\mathcal {D}^{2}_{\beta}} \gtrsim 2^{N(\alpha-\beta-1)}\sum_{n=2^N-1}^{2^{N+1}}|b_{n+1}|^{2}.$$
	
	$(2)\Rightarrow (1).$  Take  $f(z)=\sn a_{n}z^{n}\in \mathcal {D}_{\alpha}$. Set
	$$f_{+}(z)=\sn |a_{n}|z^{n}, \ \ \ z\in \dd.$$
	We have that $f_{+}\in \mathcal {D}^{2}_{\alpha}$ and $||f||_{\mathcal {D}_{\alpha}}=||f_{+}||_{\mathcal {D}^{2}_{\alpha}}$.   Now, 
	\begin{align*}
	||\hg(f)||^{2}_{ \mathcal {D}^{2}_{\beta}} & \asymp \sn (n+1)^{3-\beta}|b_{n+1}|^{2}\left|\sk \frac{a_{k}}{n+k+1} \right|^{2}\\
	& \leq \sn (n+1)^{3-\beta}|b_{n+1}|^{2}\left(\sk \frac{|a_{k}|}{n+k+1} \right)^{2}\\
	&\leq  \sn \left( \sum_{k=2^{n}-1}^{2^{n+1}}(k+1)^{3-\beta}|b_{k+1}|^{2}\left( \sum_{j=0}^{\infty}\frac{|a_{j}|}{2^{n}+j}\right) ^{2}\right) \\
	& \lesssim  \sn 2^{n(3-\beta)} \left( \sum_{j=0}^{\infty}\frac{|a_{j}|}{2^{n}+j}\right) ^{2}\left( \sum_{k=2^{n}-1}^{2^{n+1}}|b_{k+1}|^{2}\right) 
	\end{align*}
	Since $ \sum_{k=2^{n}-1}^{2^{n+1}}|b_{k+1}|^{2} =O(2^{n(\beta-\alpha-1)})$, this implies  that 
	\begin{align*}
	||\hg(f)||^{2}_{ \mathcal {D}^{2}_{\beta}} & \lesssim \sn 2^{n(2-\alpha)} \left( \sum_{j=0}^{\infty}\frac{|a_{j}|}{2^{n}+j}\right) ^{2} \\
	& \lesssim \sn 2^{n(1-\alpha)}\sum_{k=2^{n}-1}^{2^{n+1}}\left(\sum_{j=0}^{\infty} \frac{|a_{j}|}{k+j+1}\right)^{2}\\
	& \asymp \sn (n+1)^{1-\alpha}\left(\sk\frac{|a_{k}|}{n+k+1} \right)^{2}. \\
	% & \asymp ||\mathcal{H}(f_{+})||^{2}_{\mathcal {D}^{2}_{\alpha}}\lesssim ||f_{+}||^{2}_{\mathcal {D}^{2}_{\alpha}} =||f||^{2}_{\mathcal {D}_{\alpha}}.
	\end{align*}
Notice that 
	$$\mathcal{H}(f_{+})(z)=\sn \left( \sk \frac{|a_{k}|}{n+k+1}\right)z^{n}, \ \ z\in \dd. $$
	Since the Hilbert operator $\mathcal{H}$ is bounded on $\mathcal {D}^{2}_{\alpha}$ for $0<\alpha<2$ (see \cite[Theorem 1.2]{dia3}), this implies that  $$ ||\mathcal{H}(f_{+})||^{2}_{\mathcal {D}^{2}_{\alpha}} \asymp \sn (n+1)^{1-\alpha}\left(\sk\frac{|a_{k}|}{n+k+1} \right)^{2}.$$
	Therefore, we obtain 
	$$ ||\hg(f)||^{2}_{ \mathcal {D}^{2}_{\beta}} \lesssim  ||\mathcal{H}(f_{+})||^{2}_{\mathcal {D}^{2}_{\alpha}}\lesssim ||f_{+}||^{2}_{\mathcal {D}^{2}_{\alpha}} =||f||^{2}_{\mathcal {D}_{\alpha}} .$$
	The proof is complete.	\end{proof}

\begin{theorem}\label{3}
	Let  $g(z)=\sn b_{n}z^{n}\in \hd$. If $0<\alpha<2$ and $\beta \in \mathbb{R}$,
	then the following conditions are equivalent:
	
	(1) The operator $\hg:\mathcal {D}^{2}_{\alpha} \tto  \mathcal {D}^{2}_{\beta}$ is compact.
	
	(2) $\sum_{n=2^{N}}^{2^{N+1}-1}|b_{n+1}|^{2}=o\left( 2^{N(\beta-\alpha-1)}\right) $.
\end{theorem}

\begin{proof}
		$(2)\Rightarrow (1).$  
%Let
%	\begin{equation}\label{7}
%	A_{N}=\sum_{n=N}^{\infty}n^{1-\beta}|b_{n}|^{2}.
%	\end{equation}
%	Then $\{A_{N}\} \tto 0$, as $N\rightarrow \infty$. 
	For $N \in \mathbb{N}$, let  $\hg^{N}:   \mathcal {D}^{2}_{\alpha}\tto  \mathcal {D}^{2}_{\beta}$ be the operator defined as follows:
	If  $f(z)=\sn a_{n}z^{n}\in \mathcal {D}^{2}_{\alpha}$,
	\begin{equation*}
	\hg^{N}(f)(z)=\sum_{n=0}^{N}(n+1)b_{n+1}\left(\sk \frac{a_{k}}{n+k+1}\right)z^{n}.
	\end{equation*}
	The operators  $\hg^{N}$ are finite rank operators from  $\mathcal {D}^{2}_{\alpha}$ to $ \mathcal {D}^{2}_{\beta}$.  
	%Thus, the compactness of  $\hg:\mathcal {W} \tto  \mathcal {D}^{2}_{\beta}$	%will follows from

For any given $\varepsilon>0$,  there exists  $N\in\mathbb{N}$ such that 
\begin{equation}\label{s}
\sum_{k=2^{n}-1}^{2^{n+1}}|b_{n+1}|^{2} \leq \varepsilon 2^{n(\beta-\alpha-1)}, \ \ \ \mbox{for all}  \ n \geq N.		\end{equation}
	For  $f(z)=\sn a_{n}z^{n}\in \mathcal {D}^{2}_{\alpha}$, arguing as in the proof of the implication  $(2)\Rightarrow (1)$ of Theorem  \ref{th2} and using  (\ref{s}), we obtain 
	\begin{align*}
	||\hg(f)-\hg^{N}(f)||_{\mathcal {D}^{2}_{\beta}}	
	& \lesssim  \sum_{n=N}^{\infty} (n+1)^{3-\beta}|b_{n+1}|^{2}\left(\sk \frac{|a_{k}|}{n+k+1} \right)^{2}\\
	&\leq  \sum_{n=N}^{\infty}  \left( \sum_{k=2^{n}-1}^{2^{n+1}}(k+1)^{3-\beta}|b_{k+1}|^{2}\left( \sum_{j=0}^{\infty}\frac{|a_{j}|}{2^{n}+j}\right) ^{2}\right) \\
&  \lesssim   \varepsilon\sum_{n=N}^{\infty} 2^{n(1-\alpha)}\sum_{k=2^{n}-1}^{2^{n+1}}\left(\sum_{j=0}^{\infty} \frac{|a_{j}|}{k+j+1}\right)^{2}\\
&  \lesssim   \varepsilon\sum_{n=0}^{\infty} 2^{n(1-\alpha)}\sum_{k=2^{n}-1}^{2^{n+1}}\left(\sum_{j=0}^{\infty} \frac{|a_{j}|}{k+j+1}\right)^{2}\\
& \lesssim  \varepsilon ||\mathcal{H}(f_{+})||^{2}_{\mathcal {D}^{2}_{\alpha}}\lesssim \varepsilon ||f||_{\mathcal {D}^{2}_{\alpha}}.
	\end{align*}
		Hence, $\hg$ is the limit of  the  finite rank operators $\hg^{N}$  in the operator norm and therefore $\hg:  \mathcal {D}^{2}_{\alpha} \tto  \mathcal {D}^{2}_{\beta}$ is compact.
		
		$(1)\Rightarrow (2).$  As in the
		proof of Theorem  \ref{th2}, for $\frac{1}{2}<a<1$, let
			$$f_{a}(z)=(1-a)^{\frac{\alpha}{2}}\sum_{n=0}^{\infty}(n+1)^{\alpha-1}a^{n}z^{n},\ \ \ z\in \dd.$$
			Then $||f_{a}||_{ \mathcal {D}^{2}_{\alpha} }\asymp 1$ and 
			$f_{a} \tto 0$ as $a \tto 1^{-}$, uniformly
			in compact subsets of $\dd$. It follows that 
			\begin{equation}\label{L}
		\lim_{a \tto 1^-}||\hg(f_{a})||^{2}_{\mathcal {D}^{2}_{\beta}}=0.
			\end{equation}
	The proof of the implication $(1)\Rightarrow (2)$ in Theorem \ref{th2} showed  that
	$$||\hg(f_{a})||^{2}_{\mathcal {D}^{2}_{\beta}} \gtrsim (1-a)^{\alpha}2^{N(2\alpha-\beta-1)}a^{2^{N+1}}\sum_{n=2^N-1}^{2^{N+1}}|b_{n+1}|^{2}.$$
	Taking  $a=1-\frac{1}{2^{N}}$ and using \eqref{L}, we obtain that
	$$\sum_{n=2^{N}-1}^{2^{N+1}}|b_{n+1}|^{2}=o\left( 2^{N(\beta-\alpha-1)}\right) .$$
 The proof is complete. 
\end{proof}

%	\begin{corollary}\label{co1}
%		Suppose that $0<\alpha<2$ and $\beta \in\mathbb{R}$. 
%		
%	(1) If $\beta\geq\alpha+2$, then  the operator $\hg:\mathcal {D}^{2}_{\alpha} \tto  \mathcal {D}^{2}_{\beta}$ is bounded.
	
%	(2)  If $\beta<\alpha+2$, then the operator $\hg:\mathcal {D}^{2}_{\alpha} \tto  \mathcal {D}^{2}_{\beta}$ is bounded if and only if  $\sum_{n=2^{N}}^{2^{N+1}-1}|b_{n+1}|^{2}=O\left( 2^{N(\beta-\alpha-1)}\right) $.
%\end{corollary}

%As a consequence of Theorem \ref{th2}, we obtain the following result.

\begin{corollary}\label{co1}
	Suppose that   $0<\alpha<2$ and $\beta \in\mathbb{R}$. Then the operator $\mathcal{H}:\mathcal {D}^{2}_{\alpha} \tto  \mathcal {D}^{2}_{\beta}$ is bounded if and only if $\beta\geq \alpha$.
\end{corollary}

 Fot $g(z)=\sn b_{n}z^{n}\in \hd$, it is known  (see e.g., \cite[Theorem 3.1]{mat}) that $g\in \Lambda^{2}_{1/2}$ (resp. $g\in \lambda^{2}_{1/2}$) if and only if 
 $$\sum_{k=2^{n}-1}^{2^{n+1}}|b_{k}|^{2}=O(2^{-n}),\ \ \ \  (\mbox{resp}. \sum_{k=2^{n}-1}^{2^{n+1}}|b_{k}|^{2}=  o(2^{-n})).$$

Therefore, we can obtain the following corollary.
\begin{corollary}\label{co2}
	Let  $g(z)=\sn b_{n}z^{n}\in \hd$	and let  $0<\alpha<2$. Then the operator $\hg$ is bounded (resp. compact) on  $\mathcal {D}^{2}_{\alpha}$  if and only if $g\in \Lambda^{2}_{1/2}$ (resp. $g\in\lambda^{2}_{1/2}$).
\end{corollary}

 For $\alpha<0$ and $\beta \in \mathbb{R}$, the  boundedness and compactness of $\hg:\mathcal {D}^{2}_{\beta}\tto  \mathcal {D}^{2}_{\beta}$ are equivalent. As shown in the following theorem.  
 \begin{theorem}\label{th2.1}
 	Let  $g(z)=\sn b_{n}z^{n}\in \hd$.  If  $\beta \in \mathbb{R}$,
 	then the following conditions are equivalent:
 	
 	(1) The operator $\hg:\mathcal {W} \tto  \mathcal {D}^{2}_{\beta}$ is bounded.
 	
 	(2) The operator $\hg:\mathcal {W} \tto  \mathcal {D}^{2}_{\beta}$ is compact.
 	
 	(3) $\hg(1)\in \mathcal {D}^{2}_{\beta}$.
 	
 	(4) $\sum_{n=1}^{\infty} n^{1-\beta}|b_{n}|^{2}<\infty.$
 \end{theorem}
 \begin{proof}
 	$(1)\Rightarrow (3)$, $(2)\Rightarrow(1)$ and $(3)\Leftrightarrow (4)$  are obvious.
 	
 	It suffices to prove that $(4) \Rightarrow (2)$.
 	%Let $\{f_{k}\}_{k=1}^{\infty}$   be a bounded sequence in  $\mathcal {W} $  which converges to $0$ uniformly on every compact subset of  $\mathbb{D}$. Without loss of generality, we may assume   that $f_{k}(0)=0$ for all $k\geq1$ and $\sup_{k\geq 1}||f||_{\mathcal {W} }\leq 1$.
 	%To complete the proof, we have to prove that
 	%$$\lim_{k\rightarrow \infty}||\hg(f_{k})'||_{\mathcal {D}^{2}_{\beta}}=0.$$
 	%Since  $\sum_{n=1}^{\infty}n^{1-\beta}|b_{n}|^{2}<\infty$, for any $\varepsilon>0$ there exists a positive integer $N$ such that
 	For $N \in \mathbb{N}$, let
 	\begin{equation}\label{7}
 	A_{N}=\sum_{n=N}^{\infty}n^{1-\beta}|b_{n}|^{2}.
 	\end{equation}
 	%Set
 	%	$$f_{k}(z)=\sum_{j=1}^{\infty}\widehat{f_{k}}(j)z^{j},\ \ z\in \mathbb{D}.$$
 	Then $\{A_{N}\} \tto 0$, as $N\rightarrow \infty$. 
 	%For $N \in \mathbb{N}$, let  $\hg^{N}: \mathcal {W} \tto  \mathcal {D}^{2}_{\beta}$ be the operator defined as follows:
 	%If  $f(z)=\sn a_{n}z^{n}\in \mathcal {W}$,
 %	\begin{equation*}
 %	\hg^{N}(f)(z)=\sum_{n=0}^{N}(n+1)b_{n+1}\left(\sk \frac{a_{k}}{n+k+1}\right)z^{n}
 %	\end{equation*}
 %	The operators  $\hg^{N}$ are finite rank operators from  $\mathcal {W}$ to $ \mathcal {D}^{2}_{\beta}$.  %Thus, the compactness of  $\hg:\mathcal {W} \tto  \mathcal {D}^{2}_{\beta}$
 	%will follows from
 	For  $f(z)=\sn a_{n}z^{n}\in \mathcal {W}$, we have that
 	\begin{align*}
 	||\hg(f)-\hg^{N}(f)||_{\mathcal {D}^{2}_{\beta}}&\leq \sum_{n=N}^{\infty}(n+1)^{3-\beta}|b_{n+1}|^{2}\left|\sk \frac{a_{k}}{n+k+1}\right|^{2}\\
 	& \leq \sum_{n=N}^{\infty}(n+1)^{1-\beta}|b_{n+1}|^{2}\left(\sk |a_{k}|\right)^{2}\\
 	&=||f||^{2}_{\mathcal {W}} A_{N+1}  \rightarrow 0, \ \ (N\rightarrow \infty.)
 	\end{align*}
 	Hence, $\hg$ is the limit of  the  finite rank operators $\hg^{N}$  in the operator norm. Thus,   $\hg:\mathcal {W} \tto  \mathcal {D}^{2}_{\beta}$ is compact. The proof is complete.
 \end{proof}
 
 \begin{remark}
 	For $\alpha <0$, if $f(z)=\sn a_{n}z^{n}\in \mathcal {D}^{2}_{\alpha}$, then by the  Cauchy-Schwarz inequality we see that
 	\begin{equation*}
 	\sum_{n=1}^{\infty} |a_{n}|  \leq  \left(\sum_{n=1}^{\infty}|a_{n}|^{2}n^{1-\alpha}\right)^{\frac{1}{2}} \left(\sum_{n=1}^{\infty}n^{\alpha-1}\right)^{\frac{1}{2}}
 	\lesssim ||f||_{\mathcal {D}^{2}_{\alpha}}.
 	\end{equation*}
 	This implies that $\mathcal {D}^{2}_{\alpha} \subset \mathcal {W}$. In addition, for $1\leq p<\infty$, if  $f(z)=\sn a_{n}z^{n}\in S^{p}$, then  $f'(z)=
 	\sn (n+1)a_{n+1}z^{n}\in H^{p}\subset H^{1}$.
 	By Hardy's inequality we have that
 	$$\sum_{n=0}^{\infty}|a_{n+1}| \leq \pi ||f'||_{H^{1}}\lesssim ||f'||_{H^{p}}.$$
 	This shows that $S^{p}\subset \mathcal {W}$. These two facts lead us to the following two  corollaries.$\blacktriangleleft$
 	
 \end{remark}
 
 \begin{corollary}
 	Let  $g(z)=\sn b_{n}z^{n}\in \hd$.  If   $\alpha<0$ and $\beta \in \mathbb{R}$,
 	then the following statements are equivalent:
 	
 	(1) The operator $\hg: \mathcal {D}^{2}_{\alpha} \tto  \mathcal {D}^{2}_{\beta}$ is bounded.
 	
 	(2) The operator $\hg:\mathcal {D}^{2}_{\alpha} \tto  \mathcal {D}^{2}_{\beta}$ is compact.
 	
 	(3) $\sum_{n=1}^{\infty} n^{1-\beta}|b_{n}|^{2}<\infty.$
 \end{corollary}

 \begin{corollary}\label{cor2.4}
 	Let  $1\leq p<\infty$ and  let $g(z)=\sn b_{n}z^{n}\in \hd$.  The following statements are equivalent:
 	
 	(1) The operator $\hg:S^{p} \tto S^{2}$ is bounded.
 	
 	(2) The operator $\hg:S^{p} \tto S^{2}$ is compact.
 	
 	(3)$\sum_{n=1}^{\infty} n^{2}|b_{n}|^{2}<\infty.$
 \end{corollary}

 Corollary \ref{cor2.4} is also a consequence of the following theorem with $q=2$.
 
 \begin{theorem}\label{27}
 	Let  $1\leq p<\infty$ and let $g(z)=\sn b_{n}z^{n}\in \hd$.
 	
 	(1) If $1\leq q\leq 2$ and the operator  $\hg:S^{p}\rightarrow S^{q}$ is bounded, then  $\sum_{n=0}^{\infty}(n+1)^{2q-2}|b_{n}|^{q}<\infty$.

 	(2) If $2\leq q<\infty$ and $\sum_{n=0}^{\infty}(n+1)^{2q-2}|b_{n}|^{q}<\infty$, then the operator  $\hg:S^{p}\rightarrow S^{q}$ is compact.
 \end{theorem}

 To prove Theorem \ref{27}, we need some notations.
 %and fundamental lemmas.
 
 The Hardy-Littlewood space $HL(p)$ consists of those function $f\in \hd$ for which
 $$||h||^{p}_{HL(p)}=\sn (n+1)^{p-2}|\widehat{h}(n)|^{p}<\infty.$$
 It is well known (see \cite{b1}) that
 \begin{equation}\label{1}
 \pd \subset \hp \subset HL(p),\ \  0<p\leq2,
 \end{equation}
 \begin{equation}\label{2}
 HL(p) \subset \hp \subset \pd,\ \  2\leq p<\infty.
 \end{equation}

 %\begin{lemma}\label{lm2.1}
 %If $1\leq p<\infty$ and $f\in S^{p}$, then $$|f(z)|\leq \pi ||f||_{S^{p}}.$$
 %\end{lemma}

 \begin{proof of Theorem 2.9}
 	(1). Let $f(z)\equiv 1\in S^{p}$,  then $\hg(1)(z)=\sn b_{n+1}z^{n} \in S^{q}$.  The definition of $S^{p}$ shows that  $\hg(1)'\in H^{q}$. It follows from (\ref{1}) that $\hg(1)'\in H^{q}\subset HL(q)$. This implies that
 	$$\sum_{n=0}^{\infty}(n+1)^{2q-2}|b_{n}|^{q}<\infty.$$
 	
 	(2). For every $f\in S^{p}$, as in the proof of Theorem \ref{th2.1}, it  is suffices to prove that
 	\begin{equation}\label{t}
 	||\hg(f)-\hg^{N}(f)||_{S^{q}}\rightarrow 0, \ \ \ \mbox{as} \ \ \ N\rightarrow \infty.\end{equation}
 	Since $2\leq q <\infty$, by (\ref{2}) we have that $HL(q)\subset H^{q}$. To prove (\ref{t}), it suffices to prove that
 	\begin{equation*}
 	||\hg(f)'-\hg^{N}(f)'||_{HL(q)}\rightarrow 0, \ \ \ \mbox{as} \ \ \ N\rightarrow \infty.\end{equation*}
 	It is easy to check that
 	$$\hg(f)'(z)-\hg^{N}(f)'(z)=\sum_{n=N}^{\infty}(n+1)(n+2)b_{n+2}\left(\sk \frac{a_{k}}{n+k+2}\right)z^{n}.$$
 	It follows that
 	\begin{align*}
 	||\hg(f)'-\hg^{N}(f)'||^{q}_{HL(q)} &\lesssim \sum_{n=N}^{\infty} (n+1)^{3q-2}|b_{n+2}|^{q}\left|\sk \frac{a_{k}}{n+k+2}\right|^{q}\\
 	& \lesssim   \sum_{n=N}^{\infty} (n+1)^{2q-2}|b_{n+2}|^{q}\left(\sk |a_{k}|\right)^{q}\\
 	& \lesssim    ||f||^{q}_{S^{p}}\sum_{n=N}^{\infty} (n+1)^{2q-2}|b_{n+2}|^{q} \rightarrow 0 \ \ \ (N\rightarrow \infty).
 	\end{align*}
 	The proof is complete.
 \end{proof of Theorem 2.9}
 
  \begin{remark}
 	
 By mimicking the proof of Theorem \ref{th2.1}, we can easily obtain	that $\hg: \mathcal{B} \tto D^{2}_{\beta}(\beta \in \mathbb{R})$ is bounded (equivalently compact) if and only if $\sn (n+1)^{1-\beta}|b_{n}|^{2} \log^{2}(n+1)<\infty.$ 
 \end{remark}

 \begin{theorem}
 	Let  $g(z)=\sn b_{n}z^{n}\in \hd$. 
 	The following conditions are equivalent:
 	
 	(1) The operator $\hg:\mathcal {W} \tto  \mathcal {W}$ is bounded.
 	
 	(2) The operator $\hg:\mathcal {W} \tto  \mathcal {W}$ is compact.
 	
 	(3) $g\in \mathcal {W}$.
 	
 	(4) $\sum_{n=0}^{\infty} |b_{n}|<\infty.$
 \end{theorem}
 \begin{proof}
 	It suffices to prove that $(4) \Rightarrow (1)$. Arguing as in the proof of Theorem \ref{th2.1}, we have that 
 	\begin{align*}
 	||\hg(f)-\hg^{N}(f)||_{\mathcal {W}} &\leq \sum_{n=N}^{\infty} (n+1)|b_{n+1}| \left|\sk \frac{a_{k}}{n+k+1}\right|\\
 	& \leq \sum_{n=N}^{\infty} |b_{n+1}| \left(\sk |a_{k}| \right)\\
 	& =||f||_{\mathcal {W}}\sum_{n=N}^{\infty} |b_{n+1}| \rightarrow 0 \ \  (N \rightarrow \infty).
 	\end{align*}
 	The proof is complete.
 \end{proof}

 In \cite{rra}, the authors investigated the boundedness and compactness of  $\hg:\mathcal{D} \tto \mathcal{D}$. Here, we  extend the target space  $\mathcal{D}$ to  $\mathcal{D}^{2}_{\beta}$ for all $\beta \in \mathbb{R}$. Our approach is adapted from \cite{rra} with some modifications. 
 
 \begin{lemma} \label{lm2.10}\cite[Page 814]{shie} 
 	Let $f(z)=\sn a_{n}z^{n}\in \mathcal {D}$. Then there exists a positive
 	constant  $C$ independent of $f$ such that
 	$$\sum_{l=0}^{\infty}\sum_{m=0}^{\infty}\frac{|a_{l}||a_{m}|}{\log(l+m+1)}\leq C||f||^{2}_{\mathcal {D}}.$$
 	
 \end{lemma}
 
 \begin{theorem}
 	Let  $g(z)=\sn b_{n}z^{n}\in \hd$. If  $\beta \in \mathbb{R}$,
 	then the following conditions are equivalent:
 	
 	(1) The operator $\hg:\mathcal {D} \tto  \mathcal {D}^{2}_{\beta}$ is bounded.
 	
 	%(2) The operator $\hg:\mathcal {D} \tto  \mathcal {D}^{2}_{\beta}$ is compact.
 	
 	(2) $\sum_{n=N}^{\infty} n^{1-\beta}|b_{n}|^{2}=O(\frac{1}{\log N}).$
 \end{theorem}
 \begin{proof}
 	$(1)\Rightarrow(2).$  For $0<b<1$, let
 	$$f_{b}(z)=\left(\log\frac{1}{1-b}\right)^{-\frac{1}{2}}\log \frac{1}{1-bz}=\left(\log\frac{1}{1-b}\right)^{-\frac{1}{2}}\sum_{k=1}^{\infty}\frac{b^{k}}{k}z^{k}.$$
 	Then it is clear that $f_{b}\in \mathcal {D}$ for all $b\in (0,1)$ and $||f_{b}||_{\mathcal {D}}\asymp 1.$ A simple calculation shows that
 	$$\hg(f_{b})(z)=\left(\log\frac{1}{1-b}\right)^{-\frac{1}{2}}\sn(n+1)b_{n+1}\left(\sum_{k=1}^{\infty}\frac{b^{k}}{k(n+k+1)}\right)z^{n}.$$
 	Since $\hg:\mathcal {D} \tto  \mathcal {D}^{2}_{\beta}$ is bounded, so we have 
 	\begin{equation}\label{t1}
 	||\hg(f_{b})||^{2}_{\mathcal {D}^{2}_{\beta}}\gtrsim \left(\log\frac{1}{1-b}\right)^{-1}\sn (k+1)^{3-\beta}|b_{n+1}|^{2}\left(\sum_{k=1}^{\infty}\frac{b^{k}}{k(n+k+1)}\right)^{2}.
 	\end{equation}
 	For $N\in \mathbb{N}$ and $N\geq 2$, let $b_{N}=1-\frac{1}{N}$. Then it follows from (\ref{t1}) that
 	\begin{align*}
 	1&\gtrsim ||f_{b}||^{2}_{\mathcal {D}} ||\hg||^{2}_{\mathcal {D}\tto \mathcal {D}^{2}_{\beta}}\gtrsim||\hg(f_{b_{N}})||^{2}_{\mathcal {D}^{2}_{\beta}}\\
 	& \gtrsim \frac{1}{\log N}\sum_{n=N}^{\infty}(n+1)^{3-\beta}|b_{n+1}|^{2}\left(\sum_{k=1}^{N}\frac{b^{k}_{N}}{k(n+k+1)}\right)^{2}\\
 	& \gtrsim \frac{1}{\log N}\sum_{n=N}^{\infty}(n+1)^{3-\beta}|b_{n+1}|^{2} \frac{b^{N}_{N}}{(n+N+1)^{2}}\left(\sum_{k=1}^{N}\frac{1}{k}\right)^{2}\\
 	&  \gtrsim \frac{1}{\log N}\sum_{n=N}^{\infty}(n+1)^{1-\beta}|b_{n+1}|^{2}\left(\sum_{k=1}^{N}\frac{1}{k}\right)^{2}\\
 	& \gtrsim (\log N)\sum_{n=N}^{\infty}(n+1)^{1-\beta}|b_{n+1}|^{2}.
 	\end{align*}
 	This gives that 
 	$$\sum_{n=N}^{\infty}(n+1)^{1-\beta}|b_{n+1}|^{2} \lesssim \frac{1}{\log N}.$$
 	
 	$(2)\Rightarrow (1).$  Assume (2), then it is obvious that $g\in \mathcal {D}^{2}_{\beta}$.  Let $f(z)=\sn a_{n}z^{n}\in \mathcal {D}$. We  have that 
 	\begin{align*}
 	||\hg(f)||^{2}_{ \mathcal {D}^{2}_{\beta}} &\leq \sn (n+1)^{3-\beta}|b_{n+1}|^{2}\left(\sk \frac{|a_{k}|}{n+k+1}\right)^{2}\\
 	& \lesssim |a_{0}|^{2}\sn (n+1)^{1-\beta}|b_{n+1}|^{2} + \sn (n+1)^{3-\beta}|b_{n+1}|^{2}\left(\sum_{k=1}^{\infty} \frac{|a_{k}|}{n+k+1}\right)^{2}\\
 	& \lesssim ||f||^{2}_{\mathcal {D}} ||g||^{2}_{\mathcal {D}^{2}_{\beta}}+\sn (n+1)^{3-\beta}|b_{n+1}|^{2}\left(\sum_{k=1}^{\infty} \frac{|a_{k}|}{n+k+1}\right)^{2}.\\
 	\end{align*}
 	Set \begin{equation}
 	\sn (n+1)^{3-\beta}|b_{n+1}|^{2}\left(\sum_{k=1}^{\infty} \frac{|a_{k}|}{n+k+1}\right)^{2}:= S_{1}+S_{2}
 	\end{equation}
 	where 
 	$$S_{1}:=\sn (n+1)^{3-\beta}|b_{n+1}|^{2}\left(\sum_{k=1}^{n} \frac{|a_{k}|}{n+k+1}\right)^{2},$$
 	$$S_{2}:=\sn (n+1)^{3-\beta}|b_{n+1}|^{2}\left(\sum_{k=n+1}^{\infty} \frac{|a_{k}|}{n+k+1}\right)^{2}.$$
 	Next, we estimate $S_{1}$ and $S_{2}$ separately. For $S_{1}$, using the fact that $\max\{l,m\}\asymp l+m+1 (l,m \in \mathbb{N})$, we have 
 	\begin{align*}
 	S_{1} &\leq \sn (n+1)^{1-\beta}|b_{n+1}|^{2}\left(\sum_{k=1}^{n} |a_{k}|\right)^{2}\\
 	&= \sn (n+1)^{1-\beta}|b_{n+1}|^{2}\left(\sum_{l=1}^{n}\sum_{m=1}^{n}|a_{l}||a_{m}|\right)\\
 	& =\sum_{l,m=0}^{\infty}|a_{l}||a_{m}|\sum_{n=\max\{l,m\}}(n+1)^{1-\beta}|b_{n+1}|^{2}\\
 	& \lesssim  \sum_{l,m=0}^{\infty}\frac{|a_{l}||a_{m}|}{\log(l+m+1)}.
 	\end{align*}
 	By Lemma \ref{lm2.10} we have that 
 	$$S_{1}\lesssim ||f||^{2}_{\mathcal {D}}.$$
 	
 	For $S_{2}$, using the Cauchy-Schwarz inequality we obtain 
 	\begin{align*}
 	S_{2} & =\sn (n+1)^{3-\beta}|b_{n+1}|^{2}\left(\sum_{k=n+1}^{\infty} \frac{|a_{k}|k^{\frac{1}{2}}}{k^{\frac{1}{2}}(n+k+1)}\right)^{2}\\
 	& \leq  ||f||^{2}_{\mathcal {D}} \sn (n+1)^{3-\beta}|b_{n+1}|^{2} \left(\sum_{k=n+1}^{\infty} \frac{1}{k(n+k+1)^{2}}\right).
 	\end{align*}
 	Since $$\sum_{k=n+1}^{\infty} \frac{1}{k(n+k+1)^{2}}\asymp \int_{n+1}^{\infty}\frac{dx}{x(x+n+1)}.$$
 	A simple calculation shows that 
 	$$\sum_{k=n+1}^{\infty} \frac{1}{k(n+k+1)^{2}}\asymp \frac{1}{(n+1)^{2}}.$$
 	Thus, we have that 
 	$$S_{2} \lesssim ||f||^{2}_{\mathcal {D}} ||g||^{2}_{\mathcal {D}^{2}_{\beta}}.$$
 	Therefore, we deduce that 
 	$$||\hg(f)||^{2}_{ \mathcal {D}^{2}_{\beta}} \lesssim ||f||^{2}_{\mathcal {D}} .$$
 	The proof is complete.
 \end{proof}

 \begin{theorem}
 	Let  $g(z)=\sn b_{n}z^{n}\in \hd$. If  $\beta \in \mathbb{R}$,
 	then the following conditions are equivalent:
 	
 	%(1) The operator $\hg:\mathcal {D} \tto  \mathcal {D}^{2}_{\beta}$ is bounded.
 	
 	(1) The operator $\hg:\mathcal {D} \tto  \mathcal {D}^{2}_{\beta}$ is compact.
 	
 	(2) $\sum_{n=N}^{\infty} n^{1-\beta}|b_{n}|^{2}=o(\frac{1}{\log N}).$
 \end{theorem}
 \begin{proof}
  The proof is similar to Theorem 2 in \cite{rra}, with slight modifications,  we ommit the detailed.\end{proof} 
 
\iffalse \begin{theorem}
 	Suppose that $\alpha>0$ and $\beta \in \mathbb{R}$. Let  $g(z)=\sn b_{n}z^{n}\in \hd$.
 	
 	(1) If the operator $\hg:\mathcal {D}^{2}_{\alpha} \tto  \mathcal {D}^{2}_{\beta}$ is bounded, then $\sum_{n=N}^{\infty}n^{1-\beta}|b_{n}|^{2}=O(N^{-\alpha})$.
 	
 	(2) If  $g\in \mathcal {D}^{2}_{\beta-\alpha}$, then the operator $\hg:\mathcal {D}^{2}_{\alpha} \tto  \mathcal {D}^{2}_{\beta}$ is bounded.
 \end{theorem}
\fi

\section{ The range of operators $\hg$ acting on $H^{\infty}$}\label{se3}

\ \ \ \ \ \  In this section, we devote to study the range of  $\hg$ acting on $H^{\infty}$. We begin with some concepts of function spaces.

The mixed norm space $H^{p,q,\alpha}$, $0<p,q\leq \infty$, $0<\alpha<\infty$,  is the space of all
functions $f\in \hd$ for which
$$
||f||_{p,q,\alpha}=\left(\int_{0}^{1}M^{q}_{p}(r,f)(1-r)^{q\alpha-1}dr\right)^{\frac{1}{q}}<\infty, \ \mbox{for}\ 0<q<\infty,
$$
and
$$||f||_{p,\infty,\alpha}=\sup_{0\leq r<1}(1-r)^{\alpha}M_{p}(r,f)<\infty.$$

For $t\in \mathbb{R}$,  the fractional derivative of order $t$ of $f(z)=\sum_{n=0}^{\infty}a_{n}z^{n}\in \hd$ is defined by $$D^{t}f(z)=\sum_{n=0}^{\infty}(n+1)^{t}a_{n}z^{n}.$$

 If $0<p,q\leq \infty$, $0<\alpha<\infty$, then we use  $H_{t}^{p,q,\alpha}$ to denote  the space of all analytic
functions $f\in \hd$ such that
$$||D^{t}f||_{p,q,\alpha}<\infty.$$

It is a well-known fact (see \cite{pm}) that if $f\in \hd$, $0<p,q\leq \infty$, $0<\alpha,\beta<\infty$, and $s,t \in \mathbb{R}$ satisfy $s-t=\alpha-\beta$, then
$$||D^{s}f||_{p,q,\alpha} \asymp ||D^{t}f||_{p,q,\beta}.$$
Consequently, we have  $H_{s}^{p,q,\alpha}= H^{p,q,\beta}_{t}$. This fact together with 
the inclusions between mixed norm spaces (see \cite{ar}) imply that
$$H_{2}^{1,\infty,1}= H_{1+\frac{1}{p}}^{1,\infty,\frac{1}{p}}\subsetneq H^{p,\infty,1}_{1+\frac{1}{p}}=\Lambda^{p}_{\frac{1}{p}}\ \ \ \mbox{for}\ p >1.$$

Let us remark that, $H^{p,\infty,1}_{1+\frac{1}{p}}=\Lambda^{p}_{\frac{1}{p}}$ for $p>1$, and   $H^{p,\infty,1}_{1+\frac{1}{p}}=H^{1,\infty,1}_{2}$ for $p=1$. The space $\Lambda^{1}_{1}$  is not equal to the space $H^{1,\infty,1}_{2}$. As we can  easily see that  $g(z)=\log\frac{1}{1-z}\notin \Lambda^{1}_{1}$. To simplify notation, we define   $X_{p}$ as follows;
\begin{equation*}
X_{p}:=
\begin{cases}
\displaystyle{\Lambda^{p}_{\frac{1}{p}}},  &\text{ if  $p>1;$}\\
\displaystyle{H_{2}^{1,\infty,1}}, & \text{ if  $p=1.$ }\\
\end{cases}
\end{equation*}

 Concerning the action of the Hilbert operator  $\h$  on space of bounded analytic functions $H^{\infty}$,
{\L}anucha, Nowak and Pavlovic \cite{lnp} proved that the operator  $\mathcal{H}$  is bounded from $H^{\infty}$ into $BMOA$.
In fact, it is also true that
$$\mathcal {H}(H^{\infty})\subset \bigcap_{1<p<\infty}X_{p}\subset BMOA \subset \mathcal {B}.$$

In \cite{bs}, Bellavita and  Stylogiannis  investigated the exact  norm of $\h$ from  $H^{\infty}$ to $Q_{p}$ spaces, to the mean Lipschitz spaces $\Lambda^{p}_{\frac{1}{p}}$  and to certain conformally invariant Dirichlet spaces. The author of this paper also proved that the operator $\h$ is bounded from $H^{\infty}$ to the space $X_{1}$ in \cite{g1}. 
% It has the following embeddings (see [])
%$$X_{1}\subsetneq X_{p}\subsetneq BMOA \subsetneq\mathcal {B} \ \mbox{ for all}\  1<p<\infty.$$
It is easy to check that the function $g(z)=\log\frac{1}{1-z}$ belongs to $X_{1}$, and hence belongs to $X_{p}$
for each $p>1$. As mentioned previously, $\h=\hg$ with $g(z)=\log\frac{1}{1-z}$. The operator $\hg$ is a natural generalization of the   operator $\h$. In this section,
we address the question of characterizing those $g\in \hd$ for which $\hg$ is bounded form $H^{\infty}$ to the spaces  $X_{p}$ for $p\geq 1$. 
%Our main results are stated as follows.

\begin{theorem}\label{th3.1}
	Let $1\leq p<\infty$ and let $g(z)=\sn b_{n}z^{n}\in \hd$. Then the
	following two conditions are equivalent:
	
	(1) The  operator $H_{g}$ is bounded from $H^{\infty}$  to $X_{p}$.
	
	(2) $g\in X_{p}$.
\end{theorem}

To prove Theorem \ref{th3.1}, we shall use  the following elementary results.
\begin{lemma}\label{lmt}\cite[Theorem 1.2]{g1}
	The Hilbert operator $\h$ is bound from $H^{\infty}$ to $X_{1}$.
\end{lemma}

Recall that for $f(z)=\sum_{n=0}^{\infty}a_{n}z^{n}\in \hd$ and $g(z)=\sn b_{n}z^{n}\in \hd$, the Hadamard product of  $f,g\in \hd$ is defined by
$$(f\ast g)(z)=\sum_{n=0}^{\infty}a_{n}b_{n}z^{n}.$$

\begin{lemma}\cite[Lemma 3.3]{pal}\label{lm5}
	Let  $1\leq p<\infty$ and $0\leq r<1$. If $f, g\in \hd$, then
	$$M_{p}(r^{2},f\ast g)\leq M_{p}(r,f)M_{1}(r,f).$$
\end{lemma}

The following lemma provides a characterization of fractional derivatives on mixed norm spaces; see \cite[Theorem A]{blac}.
\begin{lemma}\label{lm3}
	Let $1\leq p<\infty$, $0<q\leq \infty$, $\alpha,\beta>0$ and $f\in \hd$. Then
	$$f\in H^{p,q,\alpha} \Leftrightarrow  D^{\beta}f\in H^{p,q,\alpha+\beta}.$$
\end{lemma}

\begin{proof of Theorem 3.1}
	$(1)\Rightarrow (2)$. Take $f(z)\equiv 1\in H^{\infty}$.
	Then
	$$\hg(1)(z)=\I g'(tz)dt=\frac{1}{z}\int_{0}^{z}g'(\xi)d\xi=\frac{1}{z}(g(z)-g(0)).$$
	Since $\hg(1)\in \xp$, this means  that $g(z)=z\hg(1)+g(0)\in \xp$.
	
$(2)\Rightarrow (1)$. Assume that $g\in X_{p}$. For any $f\in H^{\infty}$,  $\hg(f)$   is a
	well defined analytic function in  $\dd$ and
	\begin{align*}
	\hg(f)(z)& = \sn \left((n+1)b_{n+1}\sk \frac{a_{k}}{n+k+1}\right)z^{n}\\
	& =(g' \ast \h(f))(z).
	\end{align*}
	Now, we divided the proof into two cases.
	
	{ \bf Case $p=1$.}
	
	It obvious that
	$$D^{3}\hg(f)=(Dg' \ast D^{2}\h(f))(z).$$
	By Lemma \ref{lm5} we have that
	\begin{align}\label{eq2}
	M_{1}(r,D^{3}\hg(f))&=M_{1}(r,Dg' \ast D^{2}\h(f)) \nonumber\\
	& \lesssim M_{1}(\sqrt{r},Dg')M_{1}(\sqrt{r},D^{2}\h(f)).
	\end{align}
	It follows from Lemma \ref{lmt}  that
	\begin{equation}\label{m}
	M_{1}(\sqrt{r},D^{2}\h(f))\lesssim \frac{1}{1-r}.
	\end{equation}
	For $g(z)=\sn b_{n}z^{n}\in \hd$, it is easy to verify that
	$$\frac{g(z)-g(0)}{z}=\sn b_{n+1}z^{n}\ \ \mbox{and} \ \ Dg'(z)=\sn (n+1)^{2}b_{n+1}z^{n}.$$
	Let  $h(z)=\frac{g(z)-g(0)}{z}$,
	%\begin{equation*}
	%h(z)=
	%\begin{cases}
	%\displaystyle{\frac{g(z)-g(0)}{z}},  &\text{ if  $z\neq0;$}\\
	%\displaystyle{g'(0)}, & \text{ if $z=0.$}\\
	%\end{cases}
	%\end{equation*}
	then $h\in \hd$ and $Dg'(z)=D^{2}h(z).$ Notice that  $Dh(z)=(zh(z))'=h(z)+zh'(z)$, a simple computation shows that
	$$Dg'(z)=D^{2}h(z)=h(z)+3zh'(z)+z^{2}h''(z).$$
	By the definition of $h$, we have that
	$$h'(z)=\frac{zg'(z)-(g(z)-g(0))}{z^{2}}.$$
	and $$h''(z)=\frac{g''(z)}{z}-\frac{2g'(z)}{z^{2}}+\frac{2(g(z)-g(0))}{z^{3}}.$$
	This gives that
	$$Dg'(z)=\frac{3(g(z)-g(0))}{z}+g'(z)+zg''(z).$$
%	To estimate $M_{1}(\sqrt{r},Dg')$, it suffices to consider $r\rightarrow 1^{-}$.
	Since $g\in X_{1}$,  this implies   that
	$$M_{1}(\sqrt{r},zg'')\lesssim \frac{1}{1-r}.$$
	The inclusion $X_{1}\subsetneq \mathcal {B}$ shows that
	$$|g'(z)| \lesssim \frac{1}{1-|z|} \ \ \mbox{and} \ \  |g(z)|\lesssim \log\frac{e}{1-|z|}.$$
	It follows that
	$$M_{1}(\sqrt{r},g')\lesssim \frac{1}{1-r} \ \ \mbox{and}\ \ M_{1}(\sqrt{r},h)\lesssim \log\frac{e}{1-r}.$$
	Thus, we obtain
	$$M_{1}(\sqrt{r},Dg')\lesssim \frac{1}{1-r}.$$
	By (\ref{eq2}) and (\ref{m})  we have that
	$$M_{1}(r,D^{3}\hg(f))\lesssim \frac{1}{(1-r)^{2}}. $$
	This menas that  $D^{3}\hg(f)=D(D^{2}\hg(f))\in H^{1,\infty,2}$.
	Now,  by Lemma \ref{lm3} we have that $D^{2}\hg(f)\in H^{1,\infty,1}$. This  is equivalent to saying that  $\hg(f)\in H^{1,\infty,1}_{2}=X_{1}$.

	{ \bf Case $p>1$}
	
	It is clear that
	$$D^{2}\hg(f)=(g' \ast D^{2}\h(f))(z).$$
	For $p>1$, by Lemma \ref{lm5} we have that
	\begin{align*}
	M_{p}(r,D^{2}\hg(f))&=M_{p}(r,g' \ast D^{2}\h(f))\\
	& \lesssim M_{p}(\sqrt{r},g')M_{1}(\sqrt{r},D^{2}\h(f)).
	\end{align*}
	Since $g\in X_{p}$, so we have that
	$$ M_{p}(\sqrt{r},g')\lesssim \frac{1}{(1-r)^{1-\frac{1}{p}}}.$$
Lemma \ref{lmt} shows that
	$$M_{1}(\sqrt{r},D^{2}\h(f))\lesssim \frac{1}{1-r}.$$
	It follows that
	$$M_{p}(r,D^{2}\hg(f))\lesssim \frac{1}{(1-r)^{2-\frac{1}{p}}}.$$
	By Lemma \ref{lm3} we have 
	$g\in \Lambda^{p}_{1/p}= X_{p}$.
\end{proof of Theorem 3.1}

\begin{theorem}
	Let $g(z)=\sn b_{n}z^{n}\in \hd$. Then the
	operator $H_{g}$ is bounded from $H^{\infty}$  to $\mathcal {B}$
	if and only if $g\in \mathcal {B}$.
\end{theorem}
\begin{proof}
	The proof of the necessity  is analogous to  the Theorem \ref{th3.1}. 
	
	On the other hand. Assume that  $g\in \mathcal {B}$. For $f\in H^{\infty}$, we have 
	\begin{align*}
	\sup_{z\in \dd}(1-|z|^{2}) |\hg(f)'(z)|&=\sup_{z\in \dd}(1-|z|^{2})\left|\int_{0}^{1}tf(t)g''(tz)dt\right|\\
	& \leq ||f||_{\infty}\sup_{z\in \dd}(1-|z|^{2})\int_{0}^{1}|g''(tz)|dt\\
	& \lesssim ||f||_{\infty} ||g||_{\mathcal {B}}\sup_{z\in \dd}(1-|z|^{2}) \int_{0}^{1}\frac{dt}{(1-t|z|)^{2}}\\
	& \lesssim ||f||_{\infty} ||g||_{\mathcal {B}}.
	\end{align*}
	The proof is complete.
\end{proof}

\section{ The operators $\hg$  induced by symbols with  non-negative
	Taylor coefficients}\label{se4}

\ \ \ \ \ \ In this section, we mainly study the operators $\hg$ induced by symbols with  non-negative Taylor coefficients, acting on logarithmically weighted Bloch spaces and on Korenblum spaces.

For  $\alpha \in \mathbb{R}$,  we define the logarithmically weighted Bloch spaces $\mathcal{B}_{\log^{\alpha}}$  as follows,
$$\mathcal{B}_{\log^{\alpha}}=\left\lbrace f\in\hd : ||f||_{\mathcal{B}_{\log^{\alpha}}}=|f(0)|+\sup_{z\in \dd}(1-|z|^{2})\log^{-\alpha}\frac{e}{1-|z|^{2}} |f'(z)|<\infty\right\rbrace .$$
In particular, $\mathcal{B}_{\log^{0}}$ is just the Bolch space $\mathcal{B}$.  If $\alpha=1$, we shall denote  $\mathcal{B}_{\log^{1}}$ as $\mathcal{B}_{\log}$.

Since $\int_{0}^{1}\log\log\frac{e^{2}}{1-t}dt<\infty$ and $\int_{0}^{1}\log^{\alpha}\frac{e}{1-t}dt<\infty$  for all $\alpha \in \mathbb{R}$, this shows that $\int_{0}^{1}|f(t)|<\infty$ for every $f\in \mathcal{B}_{\log^{\alpha}}$ and hence $\hg$  is well-defined on  $\mathcal{B}_{\log^{\alpha}}$. The boundedness of $\hg$ acting between logarithmically weighted Bloch spaces is presented as follows.

\begin{theorem}\label{th4.5}
	Let $g(z)=\sn b_{n}z^{n}\in \hd$ such that $b_{n}\geq 0$ for all $n\in \mathbb{N}\cup \{0\}$. If $\alpha>-1$, then  the following statements are equivalent.
	
	(1) The operator  $\hg:\mathcal{B}_{\log^{\alpha}}\tto \mathcal{B}_{\log^{\alpha+1}}$ is bounded.
	
	(2)  $g\in \mathcal{B}$.
	
	(3) $\sum_{n=1}^{N}nb_{n}\log^{\alpha+1}(n+1)=O(N\log^{\alpha+1} N)$.
	
	(4) $\sum_{n=1}^{N}nb_{n}=O\left(N\right) $.
\end{theorem}

	To prove 	Theorem \ref{th4.5}, we  need some auxiliary lemmas. The following integral estimate was established by the author and his collaborators in \cite{zxgt}.
	
\begin{lemma}\label{s3}
Suppose that $\delta >-1$, $c>0$ and $\beta,\gamma \in \mathbb{R}$. 	If $0\leq r <1$, then
$$\int_{0}^{1}\frac{(1-t)^{\delta}}{(1-tr)^{1+\delta+c}}\log^{\beta}\frac{e}{1-t}\log^{\gamma}\frac{e}{1-tr}dt \asymp \frac{1}{(1-r)^{c}}\log^{\beta+\gamma}\frac{e}{1-r}.$$
\end{lemma}

\iffalse
\begin{lemma}
	Suppose that $\{b_{n}\}_{n=0}^{\infty}$  is a non-negative sequence. Then the following statements are equivalent:
	
	(1) $\sum_{n=1}^{N}nb_{n}\log(n+1)=O(N)$.
	
  (2)	$\sum_{n=1}^{N}nb_{n}=O\left( \frac{N}{\log N}\right) .$
\end{lemma}
\begin{proof}
  $(2)\Rightarrow (1).$  For $N\geq 2$, if $n\leq N$ then $\log(n+1)\leq \log(N+1)$.  Hence, we have that 
$$\sum_{n=1}^{N}nb_{n}\log(n+1)\leq \log(N+1) \sum_{n=1}^{N}nb_{n}\lesssim N.$$

$(1)\Rightarrow (2).$  Let $M=[\frac{N}{\log N}]$, then
$$
\sum_{n=1}^{M}nb_{n} \leq \frac{1}{\log 2}\sum_{n=1}^{M} nb_{n}\log(n+1)\\
  \lesssim M \lesssim \frac{N}{\log N}.
$$

For $M+1\leq n\leq N$,  we have that 
$$\log (n+1)\geq \log(M+1) \gtrsim \log N.$$
It follows that  
\begin{align*}
\sum_{n=M+1}^{N}a_{n}& =\sum_{n=M+1}^{N}a_{n}\log^{\alpha}(n+1)\frac{1}{\log^{\alpha}(n+1)}\\
& \leq \frac{1}{\log^{\alpha} (M+1)} \sum_{n=M+1}^{N}a_{n}\log^{\alpha}(n+1)\\
& \lesssim \frac{1}{\log^{\alpha} N}N\log^{\beta} N.
\end{align*}
Therefore, we obtain
$$\sum_{n=1}^{N}a_{n}=O(N\log^{\beta-\alpha})N.$$
The proof is complete.
\end{proof}\fi

The following lemma provides a characterization of functions in the logarithmic Bloch space with non-negative coefficients. 
\begin{lemma}\label{l4.3}\cite[Theorem 3.1]{zt}
Let  $\alpha \in \mathbb{R}$ and let $f(z)=\sn a_{n}z^{n}\in \hd$ with $a_{n}\geq 0$ for all $n\geq 0$. Then $f\in \mathcal{B}_{\log^{\alpha}}$ if and only if 
$$\sum_{n=1}^{N}na_{n}=O(N\log^{\alpha}(N+1)).$$
\end{lemma}

%\begin{lemma}\cite[Lemma 2.2]{}
%	Let  $\alpha \in \mathbb{R}$ and let  $n$ be a nonnegative integer. Then 
%	$$\int_{0}^{1}t^{n}\log^{\alpha}\frac{e}{1-t}dt \asymp \frac{\log^{\alpha}(n+1)}{n+1}.$$
%\end{lemma}

We also need the following lemma.
\begin{lemma}\label{s4}
	Let  $\{a_{n}\}_{n=0}^{\infty}$  be  a non-negative sequence. 	Suppose that $\alpha , \beta\in \mathbb{R}$. Then the following statements are equivalent:
	
	(1) $\sum_{n=1}^{N}a_{n}\log^{\alpha}(n+1)=O(N\log^{\beta}N)$.
	
	(2)	$\sum_{n=1}^{N}a_{n}=O\left(N\log^{\beta-\alpha}N \right) .$
\end{lemma}
\begin{proof} 
%	If $\alpha=0$, then there is nothing to prove.  Without loss of generality, we only need prove the case when $\alpha>0$.  If $\alpha<0$, then $-\alpha >0$. By replacing $b_{n}=a_{n}\log^{\alpha}(n+1)$ with  $a_{n}$, then the lemma remains valid.
	
	$(2)\Rightarrow (1).$  For $N\geq 2$, if  $\alpha\geq0$ and $n\leq N$, then $\log^{\alpha}(n+1)\leq \log^{\alpha}(N+1)$.  Hence, we have that 
	$$\sum_{n=1}^{N}a_{n}\log^{\alpha}(n+1)\leq \log^{\alpha}(N+1) \sum_{n=1}^{N}a_{n}\lesssim N\log^{\beta}N.$$
	
	If $\alpha<0$,  let $A_{N}=\sum_{n=1}^{N}a_{n}$, then by Abel's summation formula we have
	$$\sum_{n=1}^{N}a_{n}\log^{\alpha}(n+1)=A_{N}\log^{\alpha}(N+1)+\sum_{n=1}^{N-1}A_{n}\left(\log^{\alpha}(n+1)-\log^{\alpha}(n+2) \right) .$$
	It is easy to check that
	\begin{equation}\label{s1}
	\log^{\alpha}(n+1)-\log^{\alpha}(n+2)\asymp \frac{|\alpha|\log^{\alpha-1}(n+1)}{n}.
	\end{equation}
Since  $A_{N}=O\left( N\log^{\beta-\alpha}N\right) $,   by (\ref{s1})  we have that 
 $$\sum_{n=1}^{N-1}A_{n}\left(\log^{\alpha}(n+1)-\log^{\alpha}(n+2) \right)
\lesssim \sum_{n=1}^{N-1}\log^{\beta-1}(n+1).$$
For any  $\beta \in \mathbb{R}$, one has $$ \sum_{n=1}^{N-1}\log^{\beta-1}(n+1)=O(N\log^{\beta-1}N).$$
Therefore, we obtain 
$$\sum_{n=1}^{N}a_{n}\log^{\alpha}(n+1)=O(N\log^{\beta}N)+O(N\log^{\beta-1}N)=O(N\log^{\beta}N).$$

	$(1)\Rightarrow (2).$  If $\alpha<0$, then $\log^{\alpha}(N+1)\leq \log^{\alpha}(n+1)$ for $1\leq n\leq N$. We have that 
	$$\sum_{n=1}^{N}a_{n}\leq \log^{-\alpha}(N+1)\sum_{n=1}^{N}a_{n}\log^{\alpha}(n+1)=O\left(N\log^{\beta-\alpha}N \right).$$
For $\alpha \geq 0$,	let $M=[\frac{N}{\log^{\alpha}N}]$, then
	$$ \log M \lesssim \log\left( \frac{N}{\log^{\alpha}N} \right) \lesssim \log N$$
	and $$\log M \gtrsim \log \frac{N}{2\log^{\alpha}N}=\log N -\log (2\log^{\alpha}N)\gtrsim \log N.$$	
	This shows that $$\log M \asymp \log N.$$	
	It follows that 
	$$
	\sum_{n=1}^{M}a_{n} \leq \frac{1}{\log^{\alpha} 2}\sum_{n=1}^{M} a_{n}\log^{\alpha}(n+1)\\
	\lesssim M \log^{\beta} M
	\lesssim \frac{N}{\log^{\alpha}N} \log^{\beta}N.
	$$
	
	For $M+1\leq n\leq N$,  we have that 
	$$\log (n+1)\geq \log(M+1) \gtrsim \frac{1}{2}\log N.$$
	It follows that  
	\begin{align*}
	\sum_{n=M+1}^{N}a_{n}& =\sum_{n=M+1}^{N}a_{n}\log^{\alpha}(n+1)\frac{1}{\log^{\alpha}(n+1)}\\
	& \leq \frac{1}{\log^{\alpha} (M+1)} \sum_{n=M+1}^{N}a_{n}\log^{\alpha}(n+1)\\
	& \lesssim \frac{1}{\log^{\alpha} N}N\log^{\beta} N.
	\end{align*}
	Therefore, we obtain
	$$\sum_{n=1}^{N}a_{n}=O(N\log^{\beta-\alpha}N).$$
	The proof is complete.
\end{proof}

\begin{proof of Theorem 4.1}
$(1)\Rightarrow (3).$  Let $f(z)=\log^{\alpha+1}\frac{e}{1-z}=\sn A_{n}z^{n}$. Then by Theorem 2.31 on page 192 of the classic monograph \cite{b8}, we know that
$$A_{n} \asymp \frac{\log^{\alpha}(n+1)}{n+1}.$$
Using Lemma  \ref{l4.3} we have that  $f\in \mathcal{B}_{\log^{\alpha}}$.
Since   the operator  $\hg:\mathcal{B}_{\log^{\alpha}}\tto \mathcal{B}_{\log^{\alpha+1}}$ is bounded, this means  that
 $$\hg(f)(z)=\sn (n+1)b_{n+1}\left(\int_{0}^{1}t^{n}\log^{\alpha+1}\frac{e}{1-t}dt \right)z^{n} \in \mathcal{B}_{\log^{\alpha+1}}.$$
 Note that the coefficients of $\hg(f)$ are non-negative, by Lemma \ref{l4.3}  we  know that
 $$\sum_{n=1}^{N}(n+1)^{2}b_{n+1}\left(\int_{0}^{1}t^{n}\log^{\alpha+1}\frac{e}{1-t}dt \right)=O\left( N\log^{\alpha+1}(N+1)\right) .$$
A calculation shows that   
 $$\int_{0}^{1}t^{n}\log^{\alpha+1}\frac{e}{1-t}dt\asymp \frac{\log^{\alpha+1}(n+1)}{n+1}.$$
 Thus, we  have 
 $$\sum_{n=1}^{N}(n+1)b_{n+1}\log^{\alpha+1}(n+1)=O\left( N\log^{\alpha+1}(N+1)\right).$$
 
 $(3)\Leftrightarrow (4)$ follows from Lemma \ref{s4} and  $(2)\Leftrightarrow (4)$ follows from Lemma \ref{l4.3}. 
 
$(2)\Rightarrow (1).$ Suppose that $g\in \mathcal{B}$, then $$|g'(z)|\lesssim \frac{||g||_{\mathcal{B}}}{1-|z|^{2}}.$$
This also implies that 
\begin{equation}\label{eqg1}
|g''(z)|\lesssim \frac{||g||_{\mathcal{B}}}{(1-|z|^{2})^{2}}.
\end{equation}
For $f\in \mathcal{B}_{\log^{\alpha}}$, it is easy to check that 
$$|f(z)|\lesssim ||f||_{\mathcal{B}_{\log^{\alpha}}}\log^{\alpha+1}\frac{e}{1-|z|}.$$
By (\ref{eqg1}) and Lemma \ref{s3} we have that 
\begin{align*}
\ \ \ & ||\hg(f)||_{ \mathcal{B}_{\log^{\alpha+1}}} =|g'(0)|\left| \int_{0}^{1}f(t)dt\right| +\sup_{z\in \dd}(1-|z|^{2})\log^{-(\alpha+1)}\frac{e}{1-|z|^{2}}|\hg(f)'(z)|\\
& \lesssim  ||g||_{\mathcal{B}}||f||_{  \mathcal{B}_{\log^{\alpha}}}+||f||_{  \mathcal{B}_{\log^{\alpha}}}\sup_{z\in \dd}(1-|z|^{2})\log^{-(\alpha+1)}\frac{e}{1-|z|^{2}}\int_{0}^{1}|g''(tz)|\log^{\alpha+1}\frac{e}{1-t}dt\\
& \lesssim ||g||_{\mathcal{B}}||f||_{  \mathcal{B}}\left( 1+\sup_{z\in \dd}(1-|z|^{2})\log^{-(\alpha+1)}\frac{e}{1-|z|^{2}}\int_{0}^{1}
\frac{\log^{\alpha+1}\frac{e}{1-t}}{(1-t|z|)^{2}}dt\right) \\
& \lesssim
||g||_{\mathcal{B}_{\log}}||f||_{  \mathcal{B}}.
\end{align*}
The proof is complete.
\end{proof of Theorem 4.1}

%As mentioned previously, $\mathcal{H}_{g}(f)(z)=\int_{0}^{1}\frac{f(t)}{1-tz}dt$ with $g(z)=\log\frac{1}{1-z}$. Thus, we can obtain the following result.
\begin{corollary}\label{co4.5}
For $\alpha>-1$, the  Hilbert operator  $\mathcal{H}$ is bounded from $\mathcal{B}_{\log^{\alpha}}$ to $\mathcal{B}_{\log^{\alpha+1}}$.
\end{corollary}

Corollary \ref{co4.5}  improves and generalizes Proposition 5.2 in \cite{lnp}.
\begin{remark}
If $\alpha\leq -1$, then the Hilbert  operator  $\mathcal{H}$ is  not a bounded operator from $\mathcal{B}_{\log^{\alpha}}$ to $\mathcal{B}_{\log^{\alpha+1}}$. For $\alpha=1$, it is easy to verify that  $h(z)=\log\log\frac{e^{2}}{1-z}\in \mathcal{B}_{\log^{-1}}$. However, we have
 \begin{align*}
  \sup_{z\in  \dd}(1-|z|^{2})|\mathcal{H}(h)'(z)| & \gtrsim  \sup_{x\in [0,1)}(1-x) \int_{0}^{1}\frac{\log\log\frac{e^{2}}{1-t}}{(1-tx)^{2}}dt \\
 & \gtrsim  \sup_{x\in [0,1)}(1-x) \int_{x}^{1}\frac{\log\log\frac{e^{2}}{1-t}}{(1-tx)^{2}}dt \gtrsim \sup_{x\in [0,1)} \log\log \frac{e^{2}}{1-x} \tto \infty.
 \end{align*}
 This shows that the operator  $\mathcal{H}$ is not bounded from $\mathcal{B}_{\log^{-1}}$ to $\mathcal{B}$. For $\alpha<-1$,  take $h(z)=1 \in \mathcal{B}_{\log^{\alpha}}$,  then a similar  argument shows that $\mathcal{H}$ is not bounded from $\mathcal{B}_{\log^{\alpha}}$ to $\mathcal{B}_{\log^{\alpha+1}}$.   $\blacktriangleleft$
 
  %In [], the author proved that the Hilbert operator  $\mathcal{H}$ in series form is bounded from  $\mathcal{B}_{\log^{\alpha}}$ to $\mathcal{B}_{\log^{\alpha+1}}$ for all $\alpha \in \mathbb{R}$. 
  % Although  the integral form  and series representation of are equivalent on formal power series, their boundedness as operators when extended to spaces of analytic functions depends on the specific form of ``function magnitude'' or ``singularity'' measured by the  norm. The norm of the logarithmic Bloch space is highly sensitive to the growth rate at the boundary,  thereby amplifying the fundamental distinction between these two representations when handling certain functions.
  % This reflects the distinction between the series form and the integral form of the Hilbert operator.
\end{remark}

\begin{theorem}
	Let $g(z)=\sn b_{n}z^{n}\in \hd$ such that $b_{n}\geq 0$ for all $n\in \mathbb{N}\cup \{0\}$. Then  the following statements are equivalent:
	
	(1) The operator  $\hg:\mathcal{B}\tto \mathcal{B}$ is bounded.
	
	(2)  $g\in \mathcal{B}_{\log^{-1}}$.
	
	(3) $\sum_{n=1}^{N}nb_{n}\log(n+1)=O(N)$.
	
	(4) $\sum_{n=1}^{N}nb_{n}=O\left( \frac{N}{\log N}\right) $.
\end{theorem}

\begin{proof}
	The proof is analogous to Theorem  \ref{th4.5}, so we omit the details.\end{proof}

For $0<\alpha<\infty$, the Korenblum space  $H^{\infty}_{\alpha}$
is the space of all functions $f\in \hd$ for which
$$||f||_{H^{\infty}_{\alpha}}=\sup_{z\in \dd}(1-|z|^{2})^{\alpha}|f(z)|<\infty.$$

The Hilbert operator $\mathcal{H}$ is bounded on Korenblum space  $H^{\infty}_{\alpha}$ if and only if $0<\alpha<1$. See e.g., \cite{dj,lnp,ml}. If $0<\alpha<1$,  it turn out that the operator $\hg$ well-defined on $H^{\infty}_{\alpha}$.
\begin{proposition}
	Let $g(z)=\sn b_{n}z^{n}\in \hd$. If $0<\alpha<1$, then   the integral  $\hg(f)$ is a well defined analytic function in $\dd$ for every  $f \in H^{\infty}_{\alpha}$ and (\ref{h}) holds.
\end{proposition}
\begin{proof}
If $0<\alpha<1$, then for every $f \in H^{\infty}_{\alpha}$, $$\int^{1}_{0}|f(t)|dt\lesssim \int^{1}_{0}(1-t)^{-\alpha}dt\lesssim 1.$$
This  means that  the integral  $\hg(f)$  converges absolutely and hence $\hg(f)$ is a well defined analytic function in $\dd$.

%On the other hand, take $f_{\alpha}(z)=(1-z)^{-\alpha} \in H^{\infty}_{\alpha}$.
% If $\hg(f)$ is a well defined analytic function in $\dd$ for  every  $f \in H^{\infty}_{\alpha}$, then  $\hg(f_{\alpha})(0)$ is finite. This implies that
% $$|\hg(f_{\alpha})(0)|=|g'(0)|\int_{0}^{1}(1-t)^{-\alpha}dt<\infty
% $$
\end{proof}

The proof of the following lemma is similar to that of Lemma \ref{s4}.
\begin{lemma}\label{lm4.9}
	Let  $\{a_{n}\}_{n=0}^{\infty}$  be  a non-negative sequence. 	Suppose that  $s\geq 0$ and $t \in \mathbb{R}$. Then $\sum_{n=1}^{N}  n^{s}a_{n}=O(N^{t})$ if and only if 
	 $\sum_{n=1}^{N}  a_{n}=O(N^{t-s})$.
\end{lemma}

\begin{theorem}
	Let $g(z)=\sn b_{n}z^{n}\in \hd$ such that $b_{n}\geq 0$ for all $n\in \mathbb{N}\cup \{0\}$. Suppose that $0<\alpha<1$ and $\beta>0$. Then  the following statements are equivalent:
	
	(1) The operator  $\hg:H^{\infty}_{\alpha}\tto H^{\infty}_{\beta}$ is bounded.
	
	%(2)  $D^{\alpha}g\in H^{\infty}_{\beta}$.
	
	(2) $g\in \mathcal{B}^{1+\beta-\alpha}$.
	
	(3) $\sum_{n=1}^{N}nb_{n}=O(N^{1+\beta-\alpha})$.
	
%	(3) $\sum_{n=1}^{N}n^{1+\alpha}b_{n}=O(N^{1+\beta})$
	
\end{theorem}
\begin{proof} $(2)\Leftrightarrow (3).$ Since $b_{n}\geq 0$ for all $n\in \mathbb{N}\cup \{0\}$, it follows from Theorem  in [] that $g(z)=\sn b_{n}z^{n}\in \mathcal{B}^{\beta+1-\alpha}$ if and only if
	$$\sup_{N\geq 1}N^{\beta+1-\alpha}\sum_{n=1}^{N}nb_{n}<\infty. $$

$(1)\Rightarrow (3).$  Let $f_{\alpha}(z)=(1-z)^{-\alpha} \in H^{\infty}_{\alpha}$, then $\hg(f_{\alpha})\in H^{\infty}_{\beta}$. It is clear that
$$\hg(f_{\alpha})(z)=\sn (n+1)b_{n+1}\left(\int_{0}^{1} t^{n}(1-t)^{-\alpha}dt\right)z^{n} .$$
This gives that
$$\sn (n+1)b_{n+1}\left(\int_{0}^{1} t^{n}(1-t)^{-\alpha}dt\right)r^{n}\lesssim \frac{1}{(1-r)^{\beta}} .$$
For $N\geq 2$,  take  $r_{N}=1-\frac{1}{N}$. Then we have 
\begin{align*}
N^{\beta}&\gtrsim \sum_{n=1}^{N} (n+1)b_{n+1}\left(\int_{0}^{1} t^{n}(1-t)^{-\alpha}dt\right)r_{N}^{n}\\
& \geq \sum_{n=1}^{N} (n+1)b_{n+1}\left(\int_{0}^{1} t^{n}(1-t)^{-\alpha}dt\right)r_{N}^{N}\\
& \gtrsim  \sum_{n=1}^{N} (n+1)b_{n+1}\left(\int_{0}^{1} t^{n}(1-t)^{-\alpha}dt\right)\\
& \asymp  \sum_{n=1}^{N} (n+1)^{\alpha}b_{n+1}.
\end{align*}
This shows that 
$$\sum_{n=1}^{N} (n+1)^{\alpha}b_{n+1}=O(N^{\beta}).$$
Now,  taking  $a_{n}=(n+1)^{\alpha}b_{n}$, $t=\beta+1-\alpha$ and $s=1-\alpha$, then the desired result  follows  by  Lemma  \ref{lm4.9}.
%Note that the   coefficients of  $\hg(f_{\alpha})$ are  non-negative, by Theorem  in [] we have that
%\begin{align*}
%1 & \gtrsim \sup_{N\geq 1}N^{-(1+\beta)}\sum_{n=0}^{N}(n+1)^{2}b_{n+1}\left(\int_{0}^{1} t^{n}(1-t)^{-\alpha}dt\right)\\
%&\asymp \sup_{N\geq 1}N^{-(1+\beta)}\sum_{n=0}^{N}(n+1)^{1+\alpha}b_{n+1}.
%\end{align*} 

$(2)\Rightarrow (1).$ If $g\in \mathcal{B}^{\beta+1-\alpha}$, then we have that
$$|g'(z)|\lesssim \frac{||g||_{ \mathcal{B}^{\beta+1-\alpha}}}{(1-|z|)^{\beta+1-\alpha}}.$$
%This also implies that
%$$|g'(z)|\lesssim \frac{||g||_{ \mathcal{B}^{\beta+1-\alpha}}}{(1-|z|)^{\beta+2-\alpha}}.$$
%\begin{align*}
%||g||_{H^{\infty}_{\beta}}&=\sup_{z\in\dd}(1-|z|^{2})^{\beta}|D^{\alpha}g(z)|\\
%& \geq\sup_{x\in[0,1)}(1-x^{2})^{\beta}\sn (n+1)^{\alpha}b_{n}x^{n}.
%\end{align*}
For $f\in H^{\infty}_{\alpha}$,  by Lemma \ref{s3} we have that
\begin{align*}
||\hg(f)||_{H^{\infty}_{\beta}}& =\sup_{z\in\dd}(1-|z|^{2})^{\beta}\left|\int_{0}^{1}f(t)g'(tz)dt   \right| \\
& \leq ||f||_{H^{\infty}_{\alpha}}||g||_{ \mathcal{B}^{\beta+1-\alpha}} \sup_{z\in\dd}(1-|z|^{2})^{\beta} \int_{0}^{1}\frac{(1-t)^{-\alpha}}{(1-t|z|)^{\beta+1-\alpha}}dt\\
& \lesssim ||f||_{H^{\infty}_{\alpha}}||g||_{ \mathcal{B}^{\beta+1-\alpha}}.
\end{align*}
The proof is complete.
\end{proof}

\begin{corollary}
If $0<\alpha<1$ and $\beta>0$,  then  the following statements are equivalent:

(1)  The operator  $\mathcal{H}:H^{\infty}_{\alpha}\tto H^{\infty}_{\beta}$ is bounded.

(2) $\beta \geq \alpha$.

\end{corollary}

\section*{Conflicts of Interest}
The authors declare that there is no conflict of interest.

\section*{Funding}%% if any

The  author was supported by  the Scientific Research  Foundation of Hunan Provincial Education Department (No. 24C0222).
%\section*{Authors' Contributions}
%.
%\section*{Acknowledgements}%% if any

%\section*{Funding}%% if any
%Not applicable

%\section*{Abbreviations}%% if any
%Text for this section\ldots

\section*{Availability of data and materials}%% if any
Data sharing not applicable to this article as no datasets were generated or analysed during
the current study: the article describes entirely theoretical research.

%\section*{Ethics approval and consent to participate}%% if any
%Text for this section\ldots

%\section*{Competing interests}
%The authors declare that they have no competing interests.

%\begingroup
%\titleformat*{\section}{\fontsize{12pt}{14pt}\bfseries\selectfont}

\pdfbookmark[1]{ References}{4}
 \end{document}